\mag=\magstep1
\documentstyle{amsart}
\newtheorem{prop}{Proposition}
\newtheorem{cor}{Corollary}
\newtheorem{theo}{Theorem}

\newcommand\lra{{\longrightarrow}}

\newcommand\CC{{\Bbb C}}
\renewcommand\phi{\varphi}

\newcommand\QQ{{\Bbb Q}}
\newcommand\RR{{\Bbb R}}

\newcommand\ZZ{{\Bbb Z}}

\newcommand\cj{{\cal J}}

\newcommand\om{\omega}
\newcommand\Om{\Omega}
\newcommand\la{\lambda}

\newcommand\G{\Gamma}
\newcommand\g{\gamma}

\newcommand\k{\kappa}

\newsymbol\nmid 232D

\newcommand\de{\delta}
\newcommand\vp{\varphi}

\begin{document}
\title{A tower of genus two curves related to the Kowalewski top}
\author{Franck Leprevost}
\address{F.L.: Technische Universit\"at Berlin, Fachbereich Mathematik MA 8-1, 
Stra{\ss}e des 17. Juni 136, D-10623 Berlin Germany}
\email{leprevot@@math.tu-berlin.de}
\author{Dimitri Markushevich}
\address{D. M.: Math\'ematiques - b\^{a}t. M2, Universit\'e Lille 1,
F-59655 Villeneuve d'Ascq Cedex, France}
\email{markushe@@gat.univ-lille1.fr}

\maketitle

\abstract Several curves of genus 2 are known, such that the equations
of motion of the Kowalewski top are linearized on their Jacobians. One can expect
from transcendental approaches via solutions of equations of
motion in theta-functions, that their Jacobians are isogeneous.
The paper focuses on two such curves: Kowalewski's and that
of Bobenko--Reyman--Semenov-Tian-Shansky, the latter arising from
the solution of the problem by the method of spectral curves.
An isogeny is established between the Jacobians of these curves by purely
algebraic means, using Richelot's transformation of a genus 2 curve. 
It is shown that this isogeny respects the
Hamiltonian flows. The two curves are completed into an infinite
tower of genus 2 curves with isogeneous Jacobians.
 \endabstract




\section*{Introduction}

Several authors writing on the Kowalewski top remarked that there are a few 
apparently different curves of genus $2$ arising in the 
problem of integrating the equations of motion of the top. The one 
classically known is Kowalewski's curve \cite{Kow}; see also a modern 
exposition of her approach in \cite{A1} or \cite{A2}. The remarkable 
property of this curve $C_1$ is that the flow of solutions of the equations 
of motion is linearized on its Jacobian $J_1$, and so, the
solutions can be expressed in terms of theta-functions of two variables. 
Bobenko--Reyman--Semenov-Tian-Shansky \cite{BRS} constructed another 
curve of genus $2$ $C_2$ with the same property, but arising in a 
different way, namely, from the Lax representation for the equations of 
motion of the top. Their construction leads to a genus $2$ curve only in 
the case when the angular momentum $l$ of the top is orthogonal to the 
gravity vector $g$. So, in this case, there are
two different genus $2$ curves associated to the Kowalewski top.
It is interesting to study more closely the relation between the two
curves. \\

\noindent The authors of \cite{BRS} claim that the Jacobians 
of the two curves are isogeneous.  
They do not give an explicit proof,
but write out the solutions of the
equations of motion in terms of theta-functions on the Jacobian $J_2$
of their curve $C_2$. The formulas for the solutions define, in fact,
a map from $J_2$ onto the corresponding Liouville torus $T$. 
The knowledge of the smallest
periods of the solutions would give an information on the nature of this
map. The authors claim that this map is an isogeny. It is known also
that solutions of the equations of motion on the Jacobian
of Kowalewski's curve $J_1$ yield an isogeny from $J_1$ onto $T$. Thus,
$J_1$ and $J_2$ are isogeneous to the same abelian surface $T$, hence
isogeneous to one another. It is a natural problem to search
for an algebraic expression for such an isogeny, avoiding 
cumbersome formulas with theta-functions. An elegant and purely 
algebraic solution to this problem is given in the present paper.
\\

\noindent To our knowledge, sofar only two more curves of genus $2$, related 
to the Kowalewski top and different from the curve of Kowalewski have been 
mentioned in the literature. They are introduced 
in \cite{HvM}. The authors establish the existence of an isogeny between 
the Jacobians of these curves and of Kowalewski's. 
Their analytic approach is completely different from the one 
purely algebraic applied in the present paper, and they do not 
address the question on the relation of their curves to that of \cite{BRS}. 
\\

\noindent In the present paper, we show that the curve $C_1$ of
Kowalewski is obtained from the curve $C_2$ of
Bobenko--Reyman--Semenov-Tian-Shansky by Richelot's
transformation (\cite{R1},\cite{R2}) inducing an isogeny
of degree $4$ between their Jacobians. 
Furthermore, we
show, that in iterating Richelot's construction in a convenient way, one
can obtain a tower of countably many curves of genus $2$,
whose Jacobians are all isogeneous to that of the
curve of Kowalewski. Thus, this approach gives an infinity
of Jacobians, on which the Hamiltonian flow of the Kowalewski
top is linearized.
\\

\noindent In Section 1, we describe briefly the equations of motion of the 
Kowalewski top and the procedures leading to $C_1$ and $C_2$. We explain, why the flow of solutions of the equations 
is linearized on the Jacobians of 
the two curves. This certainly provides a linear map between the 
universal covers of the Jacobians, but still does not explain, why 
they are isogeneous.
\\

\noindent In Section 2, we describe Richelot's construction, and apply it to obtain
a (2,2)-correspondence between the curves of 
Bobenko--Reyman--Semenov-Tian-Shansky and of Kowalewski. We show
that this correspondence induces an isogeny of the Jacobians
with kernel $\ZZ /2\ZZ\oplus\ZZ /2\ZZ$, and that the isogeny
transforms the solutions of the Lax equations on $J_2$
into Kowalewski's solutions of the equations of motion of the top 
on $J_1$.
\\

\noindent In Section 3, we describe briefly Richelot's algorithm 
which leads to a tower of isogeneous abelian surfaces, and apply it 
to our situation; we obtain a tower whose ending segment
is the Jacobian of the curve of 
Bobenko et al. followed by that of the curve of Kowalewski.
\\

\noindent {\em Acknowledgements.} The second author would like
to thank A. Reyman for an excellent series of introductory 
talks at the seminar on integrable systems
in the University of Lyon in 1996/97, from which he learned
about this problem.

\section{Kowalewski's top}

\noindent We will follow the notations of \cite{BRS}. In fact, 
the integrable system introduced there is Kowa\-lewski's 
top in constant electric 
and gravitational fields, called Kowa\-lewski's
gyrostat. As soon as we are interested in the classical situation, we will
specialize all the formulas to the case when the electric field is zero.
The motion of the top can be described by the following system:
\begin{equation}\label{top}
\begin{array}{l}
\frac{dl}{dt}=[l,\om ]+[c,g],        \\
\frac{dg}{dt}=[g,\om ]                  \\
\om =\cj l
\end{array}                              \end{equation}
Here $[.,.]$ is the vector product in $\RR^3$, $l$ is the angular momentum, 
$g$ the gravity vector,
$c$ the vector of the center of mass, 
$\om$ the angular velocity,
and $\cj=I^{-1}$ the inverse of the inertia tensor 
$I=(I_{ij})_{1\leq i,j\leq 3}$, everything in a
moving frame $(e_1,e_2,e_3)$, attached to the solid. This system is Hamiltonian,
with Hamiltonian
$$
H=\frac{1}{2}(\cj l,l)-(g,c).
$$
In Kowalewski's integrable case, the inertia tensor is
$I=\mbox{diag}(1,1,1/2)$, and $c$ lies in the plane spanned by $e_1,e_2$.
If we choose the moving frame so that
the center of mass is the endpoint of $e_1$, then
$$
H=\frac{1}{2}(l_1^2+l_2^2+2l_3^2)-g_1,
$$
and there are two additional integrals of motion
$$
I_1=(l,g)^2,\;\;\; I_2=\left ({l_{{1}}}^{2}-{l_{{2}}}^{2}+2\,g_{{1}}\right )^{2}+4\,\left (l_
{{1}}l_{{2}}+g_{{2}}\right )^{2}.
$$
The equations of motion admit a Lax representation 
\begin{equation}\label{Lax.eq}\frac{dL}{dt}=
[L,M]
\end{equation}
 with Lax matrix
$$
L(\lambda )=
\left [\begin {array}{cccc} {\frac {g_{{1}}}{\lambda}}&{\frac {g_{{2}}
}{\lambda}}&-l_{{2}}+{\frac {g_{{3}}}{\lambda}}&-l_{{1}}
\\\noalign{\medskip}{\frac {g_{{2}}}{\lambda}}&-{\frac {g_{{1}}}{
\lambda}}&l_{{1}}&-l_{{2}}-{\frac {g_{{3}}}{\lambda}}
\\\noalign{\medskip}l_{{2}}+{\frac {g_{{3}}}{\lambda}}&-l_{{1}}&-2\,
\lambda-{\frac {g_{{1}}}{\lambda}}&-2\,l_{{3}}+{\frac {g_{{2}}}{
\lambda}}\\\noalign{\medskip}l_{{1}}&l_{{2}}-{\frac {g_{{3}}}{\lambda}
}&2\,l_{{3}}+{\frac {g_{{2}}}{\lambda}}&2\,\lambda+{\frac {g_{{1}}}{
\lambda}}\end {array}\right ] ,
$$
where $[L,M]=LM-ML$ for some matrix $M$, which we will not explicitize 
here. One can verify, that the invariants $H,I_1$ and $I_2$ belong to 
the algebra generated by the coefficients of $\lambda^{-2}$ and $\lambda^0$ 
in the Laurent expansions of $\mbox{Tr}(L(\la )^2)$ and
$\mbox{Tr}(L(\la )^4)$. Since these coefficients are invariant 
under the flow of  (\ref{Lax.eq}), the spectral curve
$P(\la ,\mu )=0$ is also invariant, where
$$
P(\la ,\mu )=\mbox{det}\: (L(\la )-\mu ).
$$
Let $\G$ be the non-singular compactification
of the spectral curve, and
$L(t)$ a solution of (\ref{Lax.eq}). Then we have the line
bundle $E_t$ of eigenvectors of $L(t)$ on $\G$ (it is defined a priori
on a Zariski open subset of $\G$, but it is uniquely extended to all of $\G$
as a line subbundle of a fixed vector bundle, namely,  of the trivial
one $\CC^4\times\G$). It is proved in \cite{RS}
that the evolution of the class of $E_t$ on the Jacobian of $\G$ is linear,
and the velocity $V=d[E_t]/dt$ is given by
\begin{equation}\label{velocity}
\om (V)=\sum_{p:\la (p)=\infty} \mbox{res}_p(\frac{1}{2}\mu\om)
\;\;\forall\;\;\om\in H^0(\G ,\Om^1_\G ) .
\end{equation}
Moreover, the flow is confined to the Jacobian of the curve $C_2=\G /<\tau_1>$
and parallel to the Prym variety $P(C_2/E)$, where $E=\G /<\tau_1,\tau_2>$, and
$\tau_1:(\la ,\mu )\mapsto (-\la ,\mu )\;$, $\tau_2:(\la ,\mu )\mapsto (\la ,-\mu )\;$. \\

\noindent From the physical point of view, it is natural to think of
$|g|^2$ and $I_1$ as of trivial constants of motion. The slices
$|g|^2=\g$, $I_1=\k$ represent 4-dimensional symplectic manifolds
$M_{\g\k}$ (see (1.3) of \cite{BRS} for corresponding Poisson brackets),
and the remaining first integrals $(H,I_2)$ 
yield the complete integrability of the 
Hamiltonian system on $M_{\g\k}$ in the
sense of Liouville. They define the moment map
$\mu :M_{\g\k}\lra\CC^2$, whose (compactified) fibers are
disjoint unions of Liouville tori, and the Hamiltonian flow linearizes
on their universal cover. It turns out, that the  Liouville tori 
can be identified with the Prym variety $P(C_2/E)$, if $\k\neq 0$; 
in this case, $C_2$
is of genus 3, $E$ elliptic, and $\dim P(C_2/E)=2$. If $\k =0$,
the genus of $C_2$ (resp. $E$) goes down to $2$ (resp. $0$),
and $P(C_2/E)$ becomes simply the Jacobian of $C_2$.\\

\noindent The curve $C_2$ is that of \cite{BRS} mentioned in the introduction,
and our aim is to compare it to the curve of Kowalewski.
So, we will suppose from now on 
that $I_1 =\k = 0$, and $C_2$ is of genus $2$.
We can also normalize the constants so that
$|g|^2=\g =1$. We have for $\G $ the equation
$$
\mu^4-2d_1(\la^2)\mu^2+d_2(\la^2)=0,
$$
where
$$ d_1(z)=z^{-1}-2H+2z\; ,\; d_2(z)=z^{-2}-4Hz^{-1}+I_2 ,
$$
and the equations of $C_2, E$ are obtained by substituting $\la^2=z$,
resp. $\mu^2=y$. One can check that the 1-forms
\begin{equation}\label{1-forms}
\om_0=\frac{dz}{\mu z(\mu^2-d_1(z))}\; ,\;
\om_1=\frac{1}{2}\left(\mu^2-\frac{1}{z}\right)\om_0
\end{equation}
yield a basis of $H^0(C_2 ,\Om^1_{C_2}) $. Applying (\ref{velocity}), we obtain
the following statement:
\begin{prop}
The coordinates
of the velocity vector $V$ in the basis $\om_0,\om_1$ are $(0,-1)$. 
Hence, the equation $d[E_t]/dt=V$ induces on $\mbox{\em Sym}^2(C_2)$ 
the following system:
\begin{equation}\label{Dubrovin}\begin{array}{c}
\displaystyle{\sum_{i=1,2}}\frac{dz_i/dt}{\mu_iz_i(\mu_i^2-d_1(z_i))}\; =\; 0\; , \\
\displaystyle{\sum_{i=1,2}}\frac{1}{2}\left(\mu_i^2-\frac{1}{z_i}\right)
\frac{dz_i/dt}{\mu_iz_i(\mu_i^2-d_1(z_i))}\; =\; -\: 1\; .
\end{array}\end{equation}
\end{prop}

\noindent Following \cite{BRS}, where the analogous system is 
written out for the case when $C_2$ is of genus 3, 
we will call (\ref{Dubrovin}) the Dubrovin form of the equations 
of the motion of the top.\\

\noindent The change of variables $x=\frac{1}{2}(\mu^2-z^{-1}),u=\frac{\mu}{\sqrt{2}}
(x^2+2Hx-1+\frac{1}{4}I_2)$ brings the equation of $C_2$ to the canonical
form:
$$
u^2=x(x^2+2Hx+\frac{1}{4}I_2)(x^2+2Hx-1+\frac{1}{4}I_2),
$$
and the basis (\ref{1-forms}) of $H^0(C_2 ,\Om^1_{C_2}) $ becomes
$$
\om_0=\frac{dx}{\sqrt{2}u}\;\;\; , \;\;\; \om_1=\frac{xdx}{\sqrt{2}u}\;\; .
$$
So, the Dubrovin equations can be rewritten as follows:
\begin{equation}\label{Dubrovin.2}
\frac{dx_1/dt}{u_1}+\frac{dx_2/dt}{u_2}=0\;\;\; , \;\;\;
\frac{x_1dx_1/dt}{u_1}+\frac{x_2dx_2/dt}{u_2}=-\sqrt{2}\; .
\end{equation}
The linearized equations of Kowalewski have the same form,
but on another curve of genus 2. We will write out
her solution in omitting details of calculations. 
We are using formulas from Audin \cite{A1}. Some differences
in coefficients are explained by the choice of different
dimensionless parameters: $I_{33}=1/2,|c|=|g|=1$ here, and $I_{33}=|c|=|g|=1$
in \cite{A1}. When comparing solutions, we should keep in mind that
the corresponding times are related by the equation $\tilde{t}=
\sqrt{2}t$, where $t$ is the time of Bobenko et al. (\cite{BRS}) and $\tilde{t}$ Audin's (\cite{A1}). \\

\noindent Let $x=l_1+il_2,y=l_1-il_2$. Considering them as
independent complex variables, define the new variables $\xi_1,\xi_2$ by
$$
\xi_1 = H+\frac{R(xy)-\sqrt{R(x^2)R(y^2)}}{(x-y)^2}\;\; , \;\; \xi_2 = H+\frac{R(xy)+\sqrt{R(x^2)R(y^2)}}{(x-y)^2}\;\; ,
$$

\noindent where
$$R(x)=-x^2+2Hx+1-\frac{1}{4}I_2.$$
Then the equations (\ref{top}) are reduced to the following system
\begin{equation}\label{Kow.eq}
\frac{d\xi_1/d\tilde{t}}{\eta_1}+\frac{d\xi_2/d\tilde{t}}{\eta_2}=0\;\; , \;\;
\frac{\xi_1d\xi_1/d\tilde{t}}{\eta_1}+
\frac{\xi_2d\xi_2/d\tilde{t}}{\eta_2}=i
\end{equation}
on $\mbox{Sym}^2(C_1)$, where $C_1$ is the genus $2$ curve defined by the equation
$$
\eta^2=2\xi ((\xi -H)^2+1-\frac{1}{4}I_2)((\xi -H)^2-\frac{1}{4}I_2) .
$$

\noindent Like (\ref{Dubrovin.2}), these equations describe a linearized flow
on the Jacobian of the hyperelliptic curve. Its velocity vector with
respect to the time $t$ is $V_1=(0,\sqrt{-2})$. 
It is a natural question to ask whether the two flows can be transformed
into each other by a 
holomorphic (and hence algebraic) map between the Jacobians. 
Considering the differentials
of the first kind as coordinate functions on the 
universal covering of the Jacobian of a curve,
we can represent such a map in the form $\nu_0=a\om_0, 
\nu_1=b\om_0+\sqrt{-2}\om_1$,
where $\nu_0=d\xi /\eta ,\nu_1=\xi d\xi /\eta $, 
and $a,b\in\CC$. The question is
whether it can be realized for some $a,b$ by an algebraic correspondence
between $C_1$ et $C_2$. The answer can be obtained by expressing
both solutions in terms of theta functions, but there is also a beautiful
purely algebraic construction of such a correspondence, using only 
the equations of the curves. It is described in the next section.

\section{Richelot Isogeny}

\noindent In this section, we apply Richelot's construction (\cite{R1}, 
\cite{R2}). We follow the approaches of \cite{CF}, p. 89 and of \cite{BM}. 
Let $C$ be a
genus $2$ curve defined over the ground field $K$ by an equation
$$u^2 = f(x) = G_1(x) G_2(x) G_3(x),$$
where
$$G_j(x) = g_{j2} x^2 + g_{j1}x + g_{j0} \in K[x].$$
Let $\widehat{C}$ be the genus $2$ curve defined by the following equation
$$\Delta Y^2 = F(X) = L_1(X) L_2(X) L_3(X),$$
where
$$L_1(X) = [G_2,G_3] = G'_2(X)G_3(X) - G_2(X) G'_3(X)$$
and so on, cyclically, and $\Delta = \det(g_{ij})$. A $(2,2)$-correspondence 
between $C$ and $\widehat{C}$ is defined by the curve $Z$ given over 
$C \times \widehat{C}$ by the equations

\[ \left\{ \begin{array}{ll}
   G_1(x)L_1(X) + G_2(x)L_2(X) &= 0, \\
   G_1(x)L_1(X)(x-X) &= yY.
\end{array} \right. \]
The correspondence $Z$ induces the isogeny $\varphi: J \longrightarrow
\widehat{J}$ between $J$ and $\widehat{J}$, the Jacobians of $C$ and 
$\widehat{C}$ respectively, given by the formula $\varphi([\sum n_i P_i]) = 
[\sum n_i p_2 p_1^{-1} P_i]$ for all divisor $\sum n_i P_i$ of degree zero, 
where $p_1$ (resp. $p_2$) is
the restriction to $Z$ of the projection of $C \times \widehat{C}$ to $C$ (resp.
to $\widehat{C}$). The kernel of $\varphi$ is an abelian group of type $(2,2)$,
whose non-zero elements are explicitly given in terms of the roots of the 
$G_i$'s (see \cite{BM}, p. 52). In other words, $\varphi$ is a 
$(2,2)$-isogeny of abelian surfaces, 
the so-called Richelot isogeny, and it factors the multiplication 
by $2$ on $\widehat{J}$. \\

\noindent Now let $C_1$ be Kowalewski's curve, $C=C_2$ the curve obtained by Bobenko et al., and denote by $J_1$ and $J_2$ their Jacobians respectively. These curves are defined over $K=\QQ (H,I_2)$ and, with the notations of 
the current section, the equation of $C_2$ is given by
$$u^2 = G_1(x)G_2(x)G_3(x)$$
where
$$\begin{array}{l}
G_1(x) = x, \\
G_2(x) = x^2 + 2Hx + \frac{1}{4}I_2, \\
G_3(x) = x^2 + 2Hx + \frac{1}{4}I_2-1.
\end{array}$$
\noindent It is worthwhile to compute $\Delta$ and the $L_i$'s. It follows that $
\widehat{C_2}$ is given by the equation
$$W^2=-2(X+H)[X^2 + 1 - \frac{1}{4}I_2][X^2- \frac{1}{4}I_2].$$
By the translation $\tilde{X}=X+H$, it is transformed into 
$$W^2=-2\tilde{X}[(\tilde{X}-H)^2 + 1 - 
\frac{1}{4}I_2][(\tilde{X}-H)^2- \frac{1}{4}I_2].$$
We see that $\widehat{C_2}$ is isomorphic to $C_1$ via the map
$\nu : (\tilde{X},W)\mapsto (\xi, \eta) = (\tilde{X},iW)$. So, the Jacobians $J_2$ and $J_1$ are isogeneous via the composition $\psi =\nu_\ast\circ\vp$, where $\vp$ is Richelot's isogeny, defined above. \\

\noindent There are several ways to prove that $J_1$ and $J_2$ are generically 
non-isomorphic. One of them is to compute their Igusa invariants 
\cite{Ig} and to
check that they are different. We used this procedure to complete the proof of
the following result.

\begin{theo}
$J_1$ and $J_2$ are isogeneous over $\QQ(H,I_2)(i)$ via the
isogeny $\psi$, and are generically non-isomorphic.
\end{theo}

\begin{cor}
The curves $C_1$ and $C_2$ are not isomorphic. Moreover there are no 
non-constant morphisms between $C_1$ and $C_2$.
\end{cor}

\noindent Although it directly follows from the previous theorem, the first part
of the above corollary could be proved directly. The second part is obvious, 
for a morphism between curves having the same genus $ \geq 2$ should 
be an isomorphism.

\begin{cor}
The isogeny $\psi$ transforms the flow of solutions of Dubrovin equations
(\ref{Dubrovin.2}) on $J_2$ into that of Kowalewski's equations (\ref{Kow.eq})
on $J_1$.
\end{cor}

\noindent Proof follows from the calculation of the differential
of Richelot's isogeny. Its adjoint $\de =(d_0\vp )^\ast$ can be understood
as a linear map $\de :H^0(\widehat{C_2},\Om_{\widehat{C_2}}^{1})
\lra H^0(C_2,\Om_{C_2}^1)$, and it follows from definitions that
$\de =p_{1\ast}p_2^\ast$ ($p_{1\ast}$ being the trace map for
the double covering $p_1$). This map $\de$ is computed in
\cite{BM}: $\de \left( S(X)\frac{dX}{W}\right) =S(x)\frac{dx}{u}$ for $S$ a polynomial of degree $\leq 1$.
As $d_0\nu_\ast$ is the multiplication by $-i$, we obtain:
$$
(d_0\psi )^\ast :\frac{d\xi}{\eta}\mapsto -i\frac{dx}{u}\;\; ,
\;\; (d_0\psi )^\ast :\frac{\xi d\xi}{\eta}\mapsto -i
\frac{xdx}{u}-iH\frac{dx}{u}\;\; .
$$
This implies that $d_0\psi $ transforms the generating vector
$(0,-\sqrt{2})$ of Dubrovin's flow (with respect to the basis
$\frac{dx}{u},\frac{xdx}{u}$) into $(0,\sqrt{-2})$. This ends the proof.

\section{A tower of abelian surfaces}

\noindent As explained in \cite{BM}, Richelot's method allows to construct a
tower of $(2,2)$-isogenies of abelian surfaces:

$$ \dots \longrightarrow {\cal J}_{n+1} \stackrel{\varphi_n}
{\longrightarrow} {\cal J}_n \longrightarrow \dots 
\longrightarrow {\cal J}_2 \stackrel{\varphi_1}{\longrightarrow} {\cal J}_1,$$

\noindent where ${\cal J}_n$ is the Jacobian of a genus $2$ curve ${\cal C}_n$ defined by an equation $y^2 = F_n(x)$. The algorithm takes as input a suitable factorisation of $F_{n-1}(x)= P_{n-1}(x) Q_{n-1}(x) R_{n-1}(x)$ in real polynomials of degree $2$, applies Richelot's construction on it, and outputs the polynomial $F_n(x)$ as a product $P_n(x)Q_n(x)R_n(x)$ of real polynomials of degree $2$: see \cite{BM} for a more complete description and for some applications.\\

\noindent By applying the above method with ${\cal J}_1 = J_1$ the Jacobian of Kowalewski's curve $C_1$, one obtains the following tower of isogenous Jacobians of computable curves of genus $2$

$$ \dots \longrightarrow J_{n+1} \stackrel{\psi_n}{\longrightarrow} J_n \longrightarrow \dots \longrightarrow J_2 \stackrel{\psi_1 = \psi}{\longrightarrow} J_1,$$
whose ending segment is the Jacobian of the curve of Bobenko et al. followed by the Jacobian of the curve of Kowalewski. 


\begin{thebibliography}{20}

\bibitem[1]{A1} {\em M. Audin}, Toupies.  Pr\'epublication de l'IRMA de 
Strasbourg, 1994/001.

\bibitem[2]{A2} {\em M. Audin}, Spinning tops. A course on integrable systems,
Cambridge Studies in Advanced Mathematics {\bf 51}, Cambridge Univ. Press, 1996.

\bibitem[3]{BRS}  {\em A. I. Bobenko, A. G. Reyman, and 
M. A. Semenov-Tian-Shansky}, The Kowalewski Top 99 Years Later: A Lax Pair, 
Generalizations and Explicit Solutions, Commun. Math. Phys. 
{\bf 122}  (1989), 321--354.

\bibitem[4]{BM} {\em J.-B. Bost, and J.-F. Mestre}, Moyenne 
arithm\'etico-g\'eom\'etrique et p\'eriode des courbes de genre 
$1$ et $2$, Gaz. Math. Soc. France {\bf 38} (1988), 36-64.

\bibitem[5]{CF} {\em J. W. S. Cassels, and E. V. Flynn}, Prolegomena to a 
middlebrow arithmetic of curves of genus 2, London Mathematical 
Society Lecture Note Series {\bf 230}, Cambridge Univ. Press, 1996.

\bibitem[6]{HvM} {\em E. Horozov, and P. van Moerbeke}, The Full Geometry 
of Kowalewski's Top and (1,2)-Abelian Surfaces, Commun. Pure Appl. Math. 
{\bf 42} (1989), 357--407.

\bibitem[7]{Ig} {\em J. I. Igusa}, Arithmetic variety of moduli for genus two,
Ann. of Math. {\bf 72} (1960), 612-649.

\bibitem[8]{Kow} {\em S. Kowalewski}, Sur le probl\`eme de la rotation 
d'un corps
solide autour d'un point fixe, Acta Math. {\bf 14} (1889), 177-232.

\bibitem[9]{RS}  {\em A. G. Reyman, and M. A. Semenov-Tian-Shansky}, 
Group-theoretical methods in the theory of integrable systems, In: 
Dynamical systems VII, Encycl. Math. Sci. {\bf 16} (1994), 116--225.

\bibitem[10]{R1} {\em F. Richelot}, Essai sur une m\'ethode g\'en\'erale 
pour d\'eterminer la valeur des int\'egrales ultra-elliptiques, 
fond\'ee sur des transformations remarquables de ces transcendantes,
C. R. Acad. Sci. Paris {\bf 2} (1836), 622-627.

\bibitem[11]{R2} {\em F. Richelot}, De transformatione Integralium 
Abelianorum primiordinis commentation, J. reine angew. Math. {\bf 16} (1837), 
221-341.



\end{thebibliography}
\end{document}